\documentclass[a4paper,12pt]{article}

\usepackage{amsmath}
\usepackage{amsthm}
\usepackage{amssymb}  % defines \nmid, for example
\usepackage{latexsym} % for \Box
\usepackage[all]{xy}
\usepackage{comment}
\usepackage{rotating}
\usepackage{slashbox}
\usepackage{multirow}
\usepackage{rotating}

\newcommand{\Jac}{\operatorname{Jac}}

\newcommand{\Mathematica}{{\sf Mathematica }}
\newcommand{\Sage}{{\sf Sage }}

\newcommand{\isom}{ \cong }

\newcommand{\PP}{{\mathbb P}}
\newcommand{\Q}{{\mathbb Q}}

\newcommand{\rank}{\operatorname{rank}}

\newcommand{\tor}{\operatorname{tor}}

\newcommand{\Z}{{\mathbb Z}}

\newfont{\wncyr}{wncyr10 at 12pt}
\newfont{\wncyrten}{wncyr10 at 10pt}

\newenvironment{Proof}{\par\noindent{\sc Proof:}}%
                      {\hspace*{\fill}\nobreak$\Box$\par\medskip}
                       {\hspace*{\fill}\nobreak$\Box$\par\medskip}
\newenvironment{myitemize}
{\begin{itemize}
\setlength{\itemsep}{1pt}
\setlength{\parskip}{0pt}
\setlength{\parsep}{0pt}}
{\end{itemize}}

\newtheorem{Proposition}{Proposition}[section]
\newtheorem{Theorem}[Proposition]{Theorem}
\newtheorem{Lemma}[Proposition]{Lemma}
\newtheorem{Corollary}[Proposition]{Corollary}

\theoremstyle{definition}

\newtheorem{Remark}[Proposition]{Remark}
\newtheorem{Example}[Proposition]{Example}

\newcounter{nootje}
\begin{document}
\normalsize
\title{Partitions with equal products and elliptic curves}
\author{Mohammad Sadek and Nermine El-Sissi}
\date{}
\maketitle
\let\thefootnote\relax\footnote{Mathematics Subject Classification: 14H52, 11P81}
\begin{abstract}{\footnotesize Let $a,b,c$ be distinct positive integers. Set $M=a+b+c$ and $N=abc$. We give an explicit description of the Mordell-Weil group of the elliptic curve $\displaystyle E_{(M,N)}:y^2-Mxy-Ny=x^3$ over $\Q$. In particular we determine the torsion subgroup of $E_{(M,N)}(\Q)$ and show that its rank is positive. Furthermore there are infinitely many positive integers $M$ that can be written in $n$ different ways, $n\in\{2,3\}$, as the sum of three distinct positive integers with the same product $N$ and $E_{(M,N)}(\Q)$ has rank at least $n$.}
\end{abstract}

\section{Introduction}

  The subject of partitioning integers has been used to construct infinite families of elliptic curves with positive rank. Many authors attack questions linking partitioning integers and elliptic curves arising from these integers. For example in \cite{Aguirre} the elliptic curve $E_n:y^2=x^3-nx$ where $n=a^4+b^4$ is proved to be of rank at least $2$ over $\Q(a,b)$. If $n$ can be written as a sum of two biquadrates in two different ways, then \cite{Izadi} indicates the existence of an infinite number of integers $n$ such that $E_n(\Q)$ is of rank at least $3$, and this lower bound is improved in \cite{Aguirre} to be $4$.

 In this note we study the partitions of a positive integer into three positive integers with the same product. A triple of positive integers $(x,y,z)$ is said to be a {\em partition} of a positive integer $M$ if $M=x+y+z$. The integers $x,y,z$ are the {\em parts} of the partition. We set $N=x y z$. If we are looking for all integer triples $(x,y,z)$ with sum $M$ and product $N$, then we are trying to solve two Diophantine equations. Eliminating $z$ will yield the elliptic curve \[E_{(M,N)}:y^2-Mxy-Ny=x^3\]
  In \cite{Kelly2}, the elliptic curve $E_{(M,N)}$ was shown to be of positive rank over $\Q$. We reprove this fact and give an explicit description of the torsion part $E_{(M,N)}^{\tor}(\Q)$ of the Mordell-Weil group $E_{(M,N)}(\Q)$. More precisely we show that $E_{(M,N)}^{\tor}(\Q)$ is one of the following subgroups $\Z/3\Z,\Z/6\Z$, or $\Z/2\Z\times\Z/6\Z$.

  If the curve $E_{(M,N)}$ has positive rank over $\Q$, then this means that there is an infinite number of triples of nonzero rational numbers adding up to $M$ and having product $N$. We investigate the rank of this elliptic curve if $M$ can be written as a sum of three nonzero positive integers with product $N$ in at least two different ways.

We parametrize positive integers $M$ that have two different partitions into triples with product $N$. The parametrization of $M$ is given in terms of four parameters $p,q,r,s$. Consequently we show that $E_{(M,N)}$ has rank at least $2$ over $\Q(p,q,r,s)$. Thus there exists an infinite number of pairs of integers $(M,N)$ such that $E_{(M,N)}(\Q)$ has rank at least $2$.

 A weaker result is presented when $M$ has three different partitions into triples with product $N$. An infinite parametric family of such pairs $(M,N)$ is constructed. The parametrization depends on three parameters $p,q,r$. We prove that $E_{(M,N)}$ has rank at least $3$ over $\Q(p,q,r)$, and hence the existence of an infinite number of pairs of integers $(M,N)$ for which $E_{(M,N)}(\Q)$ has rank at least $3$.

 In this note, a partitioning question is used to exhibit infinite families of elliptic curves with positive rank. This partitioning question corresponds to an interesting geometric problem.
Namely, if $(M,N)$ is a pair of positive integers such that $4M$ is the perimeter of a rectangular box $\mathcal{R}$ with integer side lengthes, and $N$ is the volume of $\mathcal{R}$, then how many different rectangular boxes with integer sides have the same perimeter and volume? Indeed, that $\rank E_{(M,N)}(\Q)\ge 1$ means that there are an infinite number of rectangular boxes with rational sides, and the same perimeter and volume. We parametrize the pairs $(M,N)$ for which there exist two rectangular boxes with integer sides, perimeter $4M$, and volume $N$. Moreover we give infinite parametric pairs $(M,N)$ for which there exists three rectangular boxes with integer sides, perimeter $4M$, and volume $N$. These pairs $(M,N)$ give rise to elliptic curves of rank at least $2$ and $3$ respectively.

\section{Partitions}

In this section we collect elementary properties about partitions with equal products.

\begin{Lemma}
\label{lem1}
Let $M$ be a positive integer that has at least two distinct partitions into triples with equal product $N$. The following statements are true:
\begin{myitemize}
\item[a)] There is no common entry between any of the triples. In particular, $N\not\in\{p,pq\}$ where $p,q$ are primes.
\item[b)] $N$ is not a prime power.
\item[c)] $N$ is a product of at least four (not necessarily distinct) primes.
\end{myitemize}
\end{Lemma}
\begin{Proof}
a) Assume that $M=a+b+c=a+d+e$ and $abc=ade$. One has $b+c=d+(bc/d)$, i.e., $bd+cd=d^2+bc$ or $b(d-c)=d(d-c)$, in other words, $b=d$ or $d=c$. In both cases, this contradicts the fact that the partitions are distinct.\\
b) If $N=p^r$, then $M=p^{r_1}+p^{r_2}+p^{r_3}=p^{s_1}+p^{s_2}+p^{s_3}$ where $r_1\ge r_2\ge r_3,s_1\ge s_2\ge s_3$. Dividing by $\min(p^{r_3},p^{s_3})$, one sum will be divisible by $p$ while the other is not. \\
c) It follows directly from a) and b).
\end{Proof}

Let $M,N,x,y,z$ be nonzero integers satisfying the following relations:
\begin{eqnarray*}
x+y+z=M \textrm{ and } xyz=N
\end{eqnarray*}
These two equations are equivalent to the following cubic equation
$Mxy-x^2y-xy^2=N$.
We homogenize the above cubic equation to obtain the following equation \[NZ^3+XY^2+X^2Y-MXYZ=0\] describing a planar curve $C_{(M,N)}$ in $\mathbb{P}_{\Q}^2$ with $(X:Y:Z)=(0:1:0)\in C_{(M,N)}(\Q)$. Therefore the Jacobian $E_{(M,N)}:=\Jac(C_{(M,N)})$ of the planar curve $C_{(M,N)}$ is an elliptic curve defined by the following Weierstrass equation. \[E_{(M,N)}:Y^2-MXY-NY=X^3\] In fact, $C_{(M,N)}$ is isomorphic to $E_{(M,N)}$ via the following transformation:
\begin{eqnarray*}
\alpha_{(M,N)}:C_{(M,N)}&\xrightarrow{\isom}&E_{(M,N)}\\
(X:Y:Z)&\mapsto&(-NZ:-NY:X)
\end{eqnarray*}
An ordered triple of nonzero integers $(d_1,d_2,d_3)$ such that $d_1+d_2+d_3=M$ and $d_1d_2d_3=N$ is sent to a point in $C_{(M,N)}(\Q)$ and hence a point in $E_{(M,N)}(\Q)$.

We will not treat triples of the form $(d,d,d)$. The reason is that the corresponding cubic curve $y^2-3dxy-d^3y=x^3$ is singular.

From now on we will assume that if $N=d^3a$, where $a$ is cube-free, then $\gcd(M,d)=1$. Otherwise the Weierstrass equation describing $E_{(M,N)}$ is not minimal. Moreover this allows us to assume that the parts of each partition of $M$ with product $N$ are coprime.

\begin{Corollary}
\label{cor1}
Let $M,N$ be nonzero integers. There is a one-to-one correspondence between the set of ordered triples
\[S_{(M,N)}(\Q)=\{(d_1,d_2,d_3):d_1+d_2+d_3=M,\;d_1d_2d_3=N,\;d_i\in\Q\}\]
and
\[E_{(M,N)}(\Q)=\{(x:y:z)\in\PP^2(\Q):y^2z-Mxyz-Nyz^2=x^3,\;xyz\ne 0\}\]
\end{Corollary}
\begin{Proof}
The bijection map is as follows:
\begin{eqnarray*}
S_{(M,N)}(\Q)&\longrightarrow& C_{(M,N)}(\Q)\longrightarrow E_{(M,N)}(\Q)\\
(d_1,d_2,d_3)&\to&(d_1:d_2:1)\to (-N:-Nd_2:d_1)\\
(-Nc/a,b/a,-a^2/bc)&\gets&(Nc:-b:-a)\gets(a:b:c)
\end{eqnarray*}
We restrict $x,y,z$ to be non-zero since we divide by them in the inverse map.
\end{Proof}

Before we proceed with investigating the torsion subgroup of the Mordell-Weil group of $E_{(M,N)}(\Q)$, we need to recall the classification of torsion points of elliptic curves over $\Q$, see (Chapter {\rm VIII}, \S 8, Theorem 7.5) in \cite{sil1}.

\begin{Lemma}
\label{lem:tors}
Let $E$ be an elliptic curve over $\Q$. Then the torsion subgroup of $E(\Q)$ is one of the following fifteen groups:
\begin{eqnarray*}
\Z/N\Z,\; 1\le N\le 10\textrm{ or }N=12;\;\Z/2\Z\times\Z/2N\Z,\textrm{ }1\le N\le 4.
\end{eqnarray*}
\end{Lemma}

\begin{Lemma}
\label{lem:EMNtors}
Let $E_{(M,N)}$ be the elliptic curve described above. One has $\Z/3\Z\subseteq E_{(M,N)}^{\tor}(\Q)$.
\end{Lemma}
\begin{Proof}
The point $(0:0:1)$ is a rational point on
\[E_{(M,N)}:y^2z-Mxyz-Nyz^2=x^3\]
Indeed the subgroup generated by $(0:0:1)$ is $\{(0:0:1),(0:N:1),O\}$.
\end{Proof}

We define the following set $S_{(M,N)}$ of classes of triples of nonzero integers as follows:
\[S_{(M,N)}=\{(a,a,b): 2a+b=M,a^2b=N\}/\sim\] where $(x_1,x_2,x_3)\sim(y_1,y_2,y_3)$ if and only if $x_i=y_j$ for some $i,j$.
\begin{Proposition}
\label{prop:tors2}
Let $a,b$ be nonzero integers such that $2a+b=M$ and $a^2b=N$. Then the point $(-N:-Na:a)$ is a torsion point of order $6$ in $E_{(M,N)}(\Q)$. In particular, if $\# S_{(M,N)}\ge 2$, then $E_{(M,N)}^{\tor}(\Q)=\Z/2\Z\times\Z/6\Z.$
\end{Proposition}
\begin{Proof}
This is direct calculation using the formulas for addition on elliptic curves, see p.58 in \cite{sil1}. More precisely, the subgroup generated by the point $(-ab,-a^2b)$ is:
\[\langle(-ab,-a^2b)\rangle=\{(-ab,-a^2b),(0,0),(-a^2,-a^3),(0,a^2b),(-ab,-ab^2),O\}\]
We observe that $(-a^2,-a^3)$ is the image of the partition $(b,a,a)$ and $(-ab,-ab^2)$ is the image of the partition $(a,b,a)$ in $E_{(M,N)}(\Q)$.
\end{Proof}

\begin{Proposition}
\label{prop:infiniterank}
Let $M,N$ be integers such that $M$ can be written as a sum of three distinct positive integers $d_1>d_2>d_3$ whose product is $N$. Assume moreover that $d_1(d_2-d_3)^3\ne d_3(d_1-d_2)^3$. Then $\rank E_{(M,N)}(\Q)\ge1$.
\end{Proposition}
\begin{Proof}
We observe that the six rational points
\begin{eqnarray*}
P_{ij}=(-d_id_j,-d_id_j^2)\in E_{(M,N)}(\Q), i\ne j,
\end{eqnarray*}
satisfy the following identities:
\begin{eqnarray*}
P_{ij}+P_{ji}=O&,&\textrm{\hskip10pt}P_{ij}+P_{ik}=(0,0),\;j\ne k,\\
P_{ij}+P_{kj}=(0,N),\;i\ne k&,&2P_{ij}=\left(\frac{d_id_j(d_i-d_k)(d_j-d_k)}{(d_i-d_j)^2},\frac{d_id_j^2(d_i-d_k)^3}{(d_i-d_j)^3}\right).
\end{eqnarray*}
We claim that $P_{ij}$ is of infinite order for every $i,j$. Assume on the contrary that $P_{ij}$ is of finite order. First of all we see that $P_{ij}$ is not a $2$-torsion, since otherwise $P_{ij}=-P_{ij}=P_{ji}$ and $d_id_j^2=d_i^2d_j$ which implies $d_i=d_j$, a contradiction. Therefore we have $12$ points of finite order $P_{ij},2P_{ij}$, $i\ne j$. Moreover, the points $P_{ij}$ and $2P_{ij}$ are distinct. For the latter statement, it is easy to show that the six points $P_{ij}$ are distinct. However, if $P_{ij}\in 2S$ where $S=\{P_{ij},1\le i,j\le 3\}$, then we have one of the following possibilities:
 if $P_{ij}=2P_{ki}$ then $P_{ki}=(0,0)$;
 if $P_{ij}=2P_{jk}$ then $P_{jk}=(0,N)$;
if $P_{ij}=2P_{ji}$ then $3P_{ij}=(0,0)$;
if $P_{ij}=2P_{ik}$ then $\displaystyle -d_j=\frac{d_k(d_i-d_j)(d_k-d_j)}{(d_i-d_k)^2}$ and observing the signs of both sides of the equality implies  $(i,j,k)\in\{(1,2,3),(3,2,1)\}$, moreover $\displaystyle -d_j^2=\frac{d_k^2(d_i-d_j)^3}{(d_i-d_k)^3}$ but this implies $(i,j,k)\in\{(2,1,3),(2,3,1)\}$, a contradiction;
if $ P_{ij}=2P_{kj}$, then $\displaystyle -d_i=\frac{d_k(d_k-d_i)(d_j-d_i)}{(d_k-d_j)^2}$ hence $(i,j,k)\in\{(2,3,1),(2,1,3)\}$, moreover $\displaystyle -d_i=\frac{d_k(d_k-d_i)^3}{(d_k-d_j)^3} $ but this implies $(i,j,k)\in\{(3,1,2),(1,3,2)\}$, a contradiction.

 Now if $2P_{ij}=2P_{ki}$ then $2(P_{ij}+P_{ik})=O$, a contradiction; if $2P_{ij}=2P_{jk}$ then $2(P_{ij}+P_{kj})=O$; if $2P_{ij}\in\{2P_{ji},2P_{ik},2P_{kj}\}$ then $\displaystyle d_1(d_2-d_3)^3=d_3(d_1-d_2)^3$ which is ruled out by our assumption. None of these $12$ points is a point of the subgroup $\{O,(0,0),(0,N)\}$. Therefore we have fifteen points of finite order, three of them make up a subgroup of order three. According to the classification of torsion points on elliptic curves, Lemma \ref{lem:tors}, this is a contradiction.
\end{Proof}

Proposition \ref{prop:infiniterank} provides an easy method to construct elliptic curves with positive rank.

\begin{Theorem}
\label{thm1}
Let $M,N$ be integers such that $M$ can be written as a sum of three positive integers $d_1>d_2>d_3$ whose product is $N$. Assume moreover that $d_1(d_2-d_3)^3\ne d_3(d_1-d_2)^3$. The Mordell-Weil group of the elliptic curve $\displaystyle E_{(M,N)}:y^2-Mxy-Ny=x^3$ satisfies $\displaystyle E_{(M,N)}(\Q)\isom \Z^r\times E_{(M,N)}^{\tor}(\Q),\;r\ge 1$, where \begin{align*}E_{(M,N)}^{\tor}(\Q)\isom\left\{\begin{array}{ll}
\Z/3\Z  & \textrm{if  $\# S_{(M,N)} =0$} \\
\Z/6\Z& \textrm{if  $\# S_{(M,N)}=1$}\\
\Z/2\Z\times\Z/6\Z & \textrm{if  $\#S_{(M,N)}=2$}
\end{array}\right.\end{align*}
\end{Theorem}
\begin{Proof}
 We recall that there is a bijection between ordered triples of nonzero integers adding up to $M$ and having product $N$, and rational points on $E_{(M,N)}$, see Corollary \ref{cor1}. Moreover Proposition \ref{prop:infiniterank} implies that if $d_1+d_2+d_3=M$ and $d_1d_2d_3=N$ where $d_1,d_2,d_3$ are distinct nonzero integers, then $(-d_id_j,-d_id_j^2)$ is a point of infinite order in $E_{(M,N)}(\Q)$. Therefore the only points of finite order are the ones coming from $S_{(M,N)}$ and the subgroup $\{O,(0,0),(0,N)\}$. Now the statement of the theorem follows from Lemma \ref{lem:EMNtors} and Proposition \ref{prop:tors2}.
\end{Proof}
\begin{Remark}
Theorem \ref{thm1} implies that no positive integer $M$ can be written in more than two different ways as $M=x_1+2x_2$ where $x_1x_2^2=N$ and $x_1,x_2>0$.
\end{Remark}

\section{A family of elliptic curves with rank at least $2$}

In Theorem \ref{thm1}, we proved that if an integer $M$ is the sum of three distinct positive integers whose product is $N$, then there is a corresponding elliptic curve $E_{(M,N)}$ of positive rank. In what follows, we study the arithmetic of these elliptic curves and introduce further conditions on the given integers to increase the rank of the corresponding elliptic curves. For this purpose, we use positive integers which have more than one partition into three distinct parts such that the product of these parts are equal in each partition.

We consider the following system of Diophantine equations:
\begin{eqnarray}
\label{eqn1}
X+Y+Z&=&U+V+W\nonumber\\
X\,Y\,Z&=&U\,V\,W
\end{eqnarray}
In \S 2 of \cite{Choudhry}, a complete solution of the following system is found.
\begin{eqnarray}
x+y+z&=&u+v+w\nonumber\\
x^3+y^3+z^3&=&u^3+v^3+w^3\label{eqn2}
\end{eqnarray}
Using the following transformations
\[\begin{array}{cc|cc}
x=-X+Y+Z, & u=-U+V+W& X=(y+z)/2,&U=(v+w)/2,\\
y=X-Y+Z,& v=U-V+W&Y=(x+z)/2,&V=(u+w)/2,\\
z=X+Y-Z,& w=U+V-W&Z=(x+y)/2,&W=(u+v)/2,
\end{array}\]
we transform the latter system (\ref{eqn2}) of equations into the following system
 \begin{eqnarray}
X+Y+Z&=&U+V+W\label{eq3}\\
(X+Y+Z)^3-24XYZ&=&(U+V+W)^3-24UVW\label{eq4}
\end{eqnarray}
Therefore, equation (\ref{eq3}) reduces equation (\ref{eq4}) to $XYZ=UVW$.

 We obtain a complete solution of system (\ref{eqn1}) using the transformations above and the complete solution of (\ref{eqn2}) found in Theorem 1 of \cite{Choudhry}. The solution is given in terms of quadratic polynomials
in four parameters such that each parameter appearing in the solution is of first degree, more explicitly, the solution of (\ref{eqn1}) is
\begin{eqnarray}\label{eq6}
X= p(r+s),& Y=q(p+s),&Z=r(q+s)\nonumber\\
U=q(r+s),&V=r(p+s),&W=p(q+s)
\end{eqnarray}
where $p,q,r,s$ are parameters.
\begin{Theorem}
\label{thm2}
Let $M$ and $N$ be positive integers given by parametrization (\ref{eq6}). Then the points $P=(-N:-NX:Y)$ and $Q=(-N:-N U:V)$ are two independent points in $E_{(M,N)}$ over $\Q(p,q,r,s)$. In particular, $E_{(M,N)}$ has rank at least $2$ over $\Q(p,q,r,s)$.
\end{Theorem}
\begin{Proof}
To show that $E_{(M,N)}$ has rank at least $2$, we need to specialize $p,q,r,s$ in the above parametrization so that the specialization of $P,Q$ are independent over $\Q$. This holds because the specialization map is a homomorphism. Putting $p=1,q=2,r=3,s=4$ yields $x=7,y=10,z=18,u=14,v=15,w=6$, and the points $P,Q$ are specialized to $(-126,-882)$ and $(-84,-1176)$ on the elliptic curve $E_{(35,1260)}$. The height-pairing matrix associated to these points has non-zero determinant $1.70464760105805$. This means that $P,Q$ are independent in $E_{(M,N)}$ over $\Q(p,q,r,s)$.
\end{Proof}

In the following example a different family of parametric pairs $(M,N)$ given in terms of one parameter $q$ is introduced such that $ E_{(M,N)}$ has rank at least $2$ over $\Q(q)$.

\begin{Example} Let $k\ge 2,q$ be positive integers. We consider the following integer $p$ given by
\[p=\frac{q^{2k-1}+1}{q+1}.\]
These integers give rise to new integers $M_k$ which can be written as follows:
\begin{eqnarray*}
M_k&=&1+pq^2+q^{2k-1}\\
&=&p+q+q^{2k}
\end{eqnarray*}
where the product of the parts in each of the above partitions of $M_k$ is $N_k=pq^{2k+1}$. The points $(-pq^{2k+1},-p^2q^{2k+3})$, $(-q^{2k+1},-q^{2k+2})$ are independent in $E_{(M_k,N_k)}$ over $\Q(q)$, and the curves $E_{(M_k,N_k)}$ have rank at least $2$ over $\Q(q)$.
\end{Example}

\section{A family of elliptic curves with rank at least $3$}

 As it has been illustrated, there is an infinite number of pairs of positive integers $(M,N)$, where $M$ can be written in two different ways as a sum of three distinct positive integers with the same product $N$, and such that for each such pair the corresponding  elliptic curve has rank at least two over $\Q$. The latter statement holds due to the fact that the specialization homomorphism is injective at infinitely many families of parameters. This suggests that if the number of partitions with the same product increases, then the rank of the corresponding elliptic curve might get larger.

We start by finding infinite number of integer solutions to the following Diophantine system:
\begin{align}
\label{eq5}
x_1+y_1+z_1=x_2+y_2+z_2=x_3+y_3+z_3\nonumber\\
x_1 y_1 z_1= x_2 y_2 z_2=x_3 y_3 z_3
\end{align}

Given positive integers $p,q,r,s$, we are looking for pairs of positive integers $(M,N)$ such that $M$ has at least three partitions into three parts with equal product $N$ and $pqrs|N$. More accurately we find $w,z$ which make the triples $(pw,qs,rz),(pqrw,s,z)$ and $(w,qrs,pz)$ solutions to the system. It is an elementary linear algebra exercise to show that
\begin{eqnarray*}
w&=&\frac{qr^2-qp-qr+p+q-r}{pqr^2-p^2qr+p^2-p-r+1}\,s,\\z&=& -\frac{p(1-qr)w+(q-1)s}{r-1}.
\end{eqnarray*}
Thus we can generate infinite number of integer solutions to the system (\ref{eq5}) by clearing denominators. More precisely the following parametrizations solve system (\ref{eq5}):
 \begin{eqnarray}
 \label{eq7}
 p&,&\;q,\;r,\nonumber\\
 s&=&pqr^2-p^2qr+p^2-p-r+1,\nonumber\\
 w&=&qr^2-qp-qr+p+q-r,\\
 z&=&pq^2r^2-pq^2r-pqr+p+q-1\nonumber.
 %&=&-(p(1-qr)(qr^2-qp-qr+p+q-r))-(q-1)(pqr^2-p^2qr+p^2-p-r+1)\\
 \end{eqnarray}

\begin{Theorem}
\label{thm3}
Let $M$ and $N$ be positive integers given by parametrization (\ref{eq7}). The points $P_1=(-prwz,-p^2rw^2z)$, $P_2=(-sz,-s^2z)$ and $P_3=(-pwz,-pw^2z)$ are three independent points in $E_{(M,N)}$ over $\Q(p,q,r)$. In particular, $ E_{(M,N)}$ has rank at least $3$ over $\Q(p,q,r)$.
\end{Theorem}
\begin{Proof}
By specializing $p=2,q=2,r=3,s=12,w=9,z=39$, we get the following partition of $159$:
\begin{align*}
159=18+24+117=108+12+39=9+72+78.
\end{align*}
The points $P_i$ are specialized to $(-2106,-37908),(-468,-5616)$ and $(-702,-6318)$ on the elliptic curve $y^2-159xy-50544y=x^3$. The determinant of the height matrix associated to these points is $4.55758994382846$. This means that $P_1,P_2,P_3$ are independent.
\end{Proof}

We can obtain parametric solutions to the system (\ref{eq5}) by solving two homogeneous linear equations, and consequently we reach an infinite family of positive integers that have three distinct partitions with the same product. These integers provide us with a family of elliptic curves of rank at least $3$.

\section{Elliptic Curves with higher rank}

The task of finding positive integers with large number of partitions into triples whose parts have equal products seems a hard problem. In fact the largest number of such partitions of a given integer up to the knowledge of the authors is $13$. The integer $17116$ has $13$ different partitions with product $2^{10} 3^3 5^27^211.13.19$, see D16 in \cite{Guy}.

We showed that positive integers that have two or three partitions with equal products can yield elliptic curves with rank at least two or three respectively. We may expect that the higher the number of partitions we can produce with equal products, the higher the rank of the corresponding elliptic curve that we can construct. Therefore if we manage to find integers with arbitrary number of such partitions, then we should predict that the corresponding elliptic curves will have arbitrary large ranks.

For example the partitions $(x_i,y_i,z_i)$ of $M=17116$ with product $N=2^{10}. 3^3.2.5^2.7^2.11.13.19$ are:

 \[\begin{array}{ccc}
     (1512,7700,7904) & (1520,7280,8316) & (1540,6840,8736) \\
     (1596,6160,9360) & (1716,5320,10080) & (1755,5120,10241) \\
     (1760,5096,10260) & (1792,4950,10374) & (2016,4180,10920) \\
     (2128,3900,11088) & (2200,3744,11172) & (2548,3168,11400)\\
     &(2736,2940,11440)&
   \end{array}
 \]
 The elliptic curve \[E:y^2-17116\,xy-2^{10} 3^3 5^2 7^2 11. 13. 19\,y=x^3\] has Mordell-Weil group $E(\Q)\isom \Z/3\Z\times \Z^6$, where the points $(-N:-Nx_i:y_i)$, $i\le 6$, are independent.

\hskip-18pt\emph{\bf{Acknowledgements.}} All calculations were performed using \Sage \cite{sage}, and \Mathematica \cite{Mathematica}.

\bibliographystyle{plain}
\footnotesize
\bibliography{partition}

\begin{thebibliography}{1}

\bibitem{Aguirre}
J.~Aguirre and J.~C. Peral.
\newblock Elliptic curves and biquadrates.
\newblock available at \texttt{http://arxiv.org/pdf/1203.2576.pdf}.

\bibitem{Choudhry}
A.~Choudhry.
\newblock Some diophantine problems concerning equal sums of integers and their
  cubes.
\newblock {\em Hardy-{R}amanujan {J}ournal}, 33:59--70, 2010.

\bibitem{Guy}
R.~K. Guy.
\newblock {\em Unsolved problems in number theory}.
\newblock Springer-Verlag, 3rd edition, 1994.

\bibitem{Izadi}
F.~A. Izadi, F.~Khoshnam, and K.~Nabardi.
\newblock Sums of two biquadrates and elliptic curves of rank $\ge 4$.
\newblock available at \texttt{http://arxiv.org/pdf/1202.5676v5.pdf}.

\bibitem{Kelly2}
John~B. Kelly.
\newblock Partitions with equal products {II}.
\newblock {\em {P}roc. {A}mer. {M}ath. {S}oc.}, 107(4):987--990, 1989.

\bibitem{Mathematica}
Wolfram Research.
\newblock {\em Mathematica Edition: Version 5.0}.
\newblock Wolfram {R}esearch.

\bibitem{sil1}
J.~Silverman.
\newblock {\em The arithmetic of elliptic curves}.
\newblock GTM 106. Springer-Verlag, New York, 1986.

\bibitem{sage}
W.\thinspace{}A. Stein et~al.
\newblock {\em {S}age {M}athematics {S}oftware ({V}ersion 4.3.5)}.
\newblock The Sage Development Team.
\newblock {\tt http://www.sagemath.org}.

\end{thebibliography}
Department of Mathematics and Actuarial Science\\ American University in Cairo\\ mmsadek@aucegypt.edu\\nelsissi@aucegypt.edu
\end{document}